\documentclass[12pt]{article}
\usepackage{amssymb}
\usepackage{geometry}
\usepackage{graphicx}
\usepackage{epstopdf}
\usepackage[english,russian]{babel}
\usepackage{amsmath,amssymb,euscript,amsthm,amsfonts,mathrsfs,amscd}
 \usepackage{color}
 \usepackage{floatflt}
 \usepackage{mathrsfs}

 \title{
 \bf On Polynomial Bounds  of Convergence for the Availability Factor }
 
\author{Alexander Veretennikov,\footnote{University of Leeds, Leeds, UK, \& National Research University Higher School of Economics,
\& Institute for Information Transmission Problems, Moscow, Russia}%
\and Galina Zverkina\footnote{Moscow State University of Railway Engineering, Moscow, Russia} \, \thanks{Both authors are supported by the RFBR, project No 14-01-00319 A.
For the first author the article was prepared within the framework of a subsidy granted to the HSE by the Government of the Russian Federation for the implementation of the Global Competitiveness Program.}}

\begin{document}

\maketitle

\begin{abstract}
A computable estimate of the readiness coefficient for a standard  binary-state system is established in the case where both working and repair time distributions possess heavy tails.

\end{abstract}

\section{Introduction}

Let us consider a restorable system, which may be either in the working state during a random time \(\xi\) with a distribution function \(F_1(s)\stackrel{\mbox{def}}{=\!\!=}
\mathbf{P}\{\xi \leqslant s\}\;\), or it may be broken down and being restored by some service during another random time \(\eta\) with a distribution function \(F_2(s)\stackrel{\mbox{def}}{=\!\!=}\mathbf{P}\{\eta \leqslant s\}\;\).
All periods of working and repairing are alternate and independent.
The readiness coefficient (or availability factor) \(A(t)\) is defined as the probability that at time \(t\) the system is in the working  (= serviceable) state.

Often in the literature it is accepted that at initial time $t=0$ the system is serviceable and that it is in the beginning of its working period. We consider a more general case assuming that the activity of the system may have started  earlier so that at $t=0$ the system can be in one of the two states: perfect functionality or complete failure; and further that {\it before  $t=0$} the system already spent time $x$ in its current state.

Let us formalize the definition of readiness coefficient (availability factor).
We assume that $ \xi_i$ are random variables with a common distribution function $F_1(x)= \mathbf{P} \{\xi_i\leqslant x \}\;$; likewise,  $ \eta_i$ are random variables with a (another) common distribution function $F_2(x)= \mathbf{P} \{\eta_i\leqslant x \}\;$; all of them are  mutually  independent.

If at time $t=0$ our system is working and its elapsed working time  before $t=0$ equals $\;x\;$, then the {\em residual time} of this working period is a random variable denoted by $\xi^{(x)}\;$; its distribution function is denoted by
$$
F_1^{(x)}(s)\stackrel{\mbox{\small def}}{=\!\!=}\mathbf{P}\{\xi^{(x)}\leqslant s\}=\mathbf{P}\{\xi \leqslant s+x|\xi>x\}=1-\frac{1-F_1(x+s)}{1-F_1(x)}\;.
$$
Correspondingly, if at time $t=0$ the system is under repair and the duration of this repair before $t=0$ equals $x\;$, then the residual time of this repair period is a random variable denoted by  $\eta^{(x)}$ with a distribution function
$$
F_2^{(x)}(s)\stackrel{\mbox{\small def}}{=\!\!=}\mathbf{P}\{\eta^{(x)}\leqslant s\}=\mathbf{P}\{\eta \leqslant s+x|\eta>x\}=1-\frac{1-F_2(x+s)}{1-F_2(x)}\;.
$$

In the first case we will use notations
$$
t_0\stackrel{\mbox{\small def}}{=\!\!=}0\;, \quad t_1 \stackrel{\mbox{\small def}}{=\!\!=} \xi^{(x)}+\eta_1\;,\quad
t_i\stackrel{\mbox{\small def}}{=\!\!=} \xi^{(x)}+\eta_1+\sum\limits_{j=2}^{i}(\xi_j+ \eta_j)
\; ,
$$
and $
t_0'\stackrel{\mbox{\small def}}{=\!\!=}\xi^{(x)}\;,\quad t_i'\stackrel{\mbox{\small def}}{=\!\!=}\xi^{(x)}+ \sum\limits_{j=1}^{i-1}(\eta_j+ \xi_{j+1})\;.
$

In the second case
$
t_1\stackrel{\mbox{\small def}}{=\!\!=}\eta^{(x)}\;, \quad t_i\stackrel{\mbox{\small def}}{=\!\!=} \eta^{(x)}+\sum\limits_{j=2}^{i}(\xi_{j-1}+ \eta_j)\; ,
$
and
$$
t_1'\stackrel{\mbox{\small def}}{=\!\!=}0\;,\quad t_1'\stackrel{\mbox{\small def}}{=\!\!=}\eta^{(x)}+\xi_1\;,\quad  t_i'\stackrel{\mbox{\small def}}{=\!\!=}\eta^{(x)}+\xi_1+ \sum\limits_{j=2}^i(\eta_j+ \xi_{j})\;.
$$

In this notation $
A (t)\stackrel{\mbox{\small def}}{=\!\!=} \mathbf{P} \left\{t \in\bigcup\limits_{i}[t_i,t_i') \right\}\;
$.

It is well known that if distributions of $\;\xi+\eta\;$ are non-arithmetical and $\mathbf{E}\,\xi+\mathbf{E}\,\eta<\infty\;$, then there exists a limiting value $$
\lim\limits _{t\to\infty}A (t)\stackrel{\mbox{\small def}}{=\!\!=} A =\frac{\mathbf{E}\,\xi} {\mathbf{E}\,\xi+\mathbf{E}\,\eta}\;.
$$
Moreover, if $ \mathbf{E}\,\xi^n+\mathbf{E}\,\eta^n<\infty\;$ for some $n>1\;$, then
$$
\limsup\limits _{t\to\infty}\big|A (t)-A \big|t^{n-1}<\infty
$$
(see e.g.  \cite[Theorem 3, Appendix 1]{VeZv1} or  \cite[Theorem 10.7.4]{VeZv2}).

In other words, for any $\alpha\in(0,n-1]$ there exists a constant  $C(\alpha)$ such that for all $t>0$
$$
\big|A (t)-A \big|\leqslant \frac{C(\alpha)}{(1+t)^\alpha}\;.
$$

However the general theory does not provide neither the value of $C(\alpha)$, nor any bound for it.

Any knowledge of the value $C(\alpha)$ or its bound is rather important in applications.
Also, in the case where nothing was known earlier about such a constant at all, even rough estimates could be useful.
 The goal of this paper is to give explicit  estimates to this constant.

This paper is an extended version of the conference publication \cite{VeZv2015}. The section 2 contains assumptions and notations; the section 3 presents the main result; the last section 4 provides the  full proof.

\section{Assumptions and notations}
\subsection{Assumptions} We suppose that $\;F_1(x)=1-\displaystyle e^{-\int\limits_0^x\lambda(s)\,\mathrm{d} x} \;$, i.e., almost everywhere\linebreak $\lambda(s)= \displaystyle \frac{F_1'(s)}{1-F_1(s)}\;$, and for some  $\Lambda>K_1>3\;$, $K_2>3\;$
\begin{eqnarray}
    \Lambda\geqslant\lambda(s)\geqslant \displaystyle\frac{K_1}{1+s}\; \mbox{ when $s>0$}\;,\label{usl1}
\\
F_2(s)\geqslant1-\displaystyle \frac{1}{(1+s)^{K_2\!\!\!\phantom{{}^{1}}}}\; \mbox{  when $s>0\;$};\label{usl3}
\end{eqnarray}
{\it we do not assume  continuity of $F_2(s)\;$}.

Note that from (\ref{usl1}) it follows that $F_1(s)\geqslant1-\displaystyle \frac{1}{(1+s)^{K_1{\phantom{^1}}}}\;$ and there is \linebreak $f_1'(s)\stackrel{\mbox{\small def}}{=\!\!=} F_1'(s)\in\left(\displaystyle \frac{K_1 e^{-\Lambda s}}{1+s}, \displaystyle \frac{\Lambda}{(1+s)^{K_1{\phantom{^1}}}\!\!}\right)\;$.

So, (\ref{usl1}) and (\ref{usl3}) imply that for all $a\in[1,K_1-1)$ and $b\in[1, K_2-1)$ we have,
\begin{eqnarray}
  \mathbf{E}\,\xi^a\geqslant \int\limits_0^\infty s^a \frac{K_1 e^{-\Lambda s}}{1+s}\,\mathrm{d}\,s \geqslant K_1 \int\limits_0^\infty s^{a-1}e^{-\Lambda s}\,\mathrm{d}\,s=\frac{K_1\Gamma(a)}{\Lambda^a}\stackrel{\mbox{\small def}}{=\!\!=} \mu_a; \label{ocsnizu}
  \\
  \mathbf{E}\,\xi^a= a\int\limits_0^\infty s^{a-1}(1-F_1(s))\,\mathrm{d}\,s\leqslant a\int\limits_0^\infty \frac{s^{a-1}}{(1+s)^{K_1{\phantom{^1}}}\!\!}\,\mathrm{d}\,s <\frac{a}{K_1-a}\stackrel{\mbox{\small def}}{=\!\!=} m_1(a)\;,\label{oc1sverxu}
  \\
  \mbox{similarly, }\mathbf{E}\,\eta^b<\frac{b}{K_2-b}\stackrel{\mbox{\small def}}{=\!\!=} m_2(b)\;\label{oc2sverxu},
\end{eqnarray}
which suffices for the existence of $A<\infty\;$.

Notice that $\lambda(s)$ is called  \emph{intensity} of failure of the recoverable system, of course, while it is working.
\subsection{Notations}\label{def}

\noindent\textbf{1. }
Denote $K\stackrel{\mbox{\small def}}{=\!\!=}\min(K_1,K_2)\;$.

\noindent \textbf{2. }The behaviour of the system under consideration may be presented by the random process
$$
X_t=(n_t,x_t)=
\left\{
\begin{array}{cc}
               (1,t-t_i)\;, & \mbox{ if } t\in[t_i,t_i')\;; \\
               \\
               (2,t-t_i')\;, &  \mbox{ if } t\in[t_i',t_{i+1})\;;
             \end{array}
\right.
\;\;\;\;
n(X_t)\stackrel{\mbox{\small def}}{=\!\!=} n_t\;,\; x(X_t)\stackrel{\mbox{\small def}}{=\!\!=} x_t\;.
$$

The state space of the process $X_t$ is a set $\mathscr{X}\stackrel{\mbox{\small def}}{=\!\!=}\{\{ 1,2\}\times\mathbb{R}_+\}$ with a standard $\sigma$-algebra.

Denote $\mathscr{S}_j \stackrel{\mbox{\small def}}{=\!\!=}\{(j,x),\,x\in \mathbb{R}_+\}\subset\mathscr{X}$ \, ($j=(1,2)$)\,\,.

Let $X_0=(n_0,x_0)\;$.

\noindent \textbf{3. }Denote (here $j=(1,2)$):
\begin{equation}
    \begin{array}{ll}
\mathrm{M}_j(k)\stackrel{\mbox{\small def}}{=\!\!=} k\displaystyle\int\limits_0^\infty\frac{s^{\,k-1}}{(1+s)^{K_j} \phantom{^{1^1}\!\!\!\!\!}}\,\mathrm{d} s\;; &
\mathrm{M}_j^{(x)}(k)\stackrel{\mbox{\small def}}{=\!\!=} \frac{k}{1-F_{j\,}(x)}\displaystyle \int\limits_0^\infty\frac{s^{\,k-1}}{(1+s+x)^{K_j} \phantom{^{1^1}\!\!\!\!\!\!}}\,\mathrm{d} s\;;
\\
 \varkappa(T)\stackrel{\mbox{\small def}}{=\!\!=}\displaystyle \int\limits_0^\infty\frac{K_1e^{-\Lambda s}}{1+T+s}\,\mathrm{d} s\;; & F_j^{(a)}(s)\stackrel{\mbox{\small def}}{=\!\!=}1-\displaystyle\frac{1-F_j(s+a)}{1-F_j(a)}\;;
 \\
 f_j^{(a)}(s)\stackrel{\mbox{def}}{=\!\!=} \left(F^{(a)}_j(s)\right)'\;;&
 \varphi_{x,y}(s)\stackrel{\mbox{def}}{=\!\!=} \min\left(f_1^{(x)}(s),f_1^{(y)}(s) \right)\;;
 \\
 \varkappa_{x,y}\stackrel{\mbox{def}}{=\!\!=} \displaystyle\int\limits_0^\infty \varphi_{x,y}(s)\,\mathrm{d} \,s
 \;; & \Phi_{x,y}(s)\stackrel{\mbox{def}}{=\!\!=} \displaystyle\int\limits_0^s\varphi_{x,y}(u)\,\mathrm{d} u\;;
 \\
  \widehat{\Phi}_{x,y}(s)\stackrel{\mbox{def}}{=\!\!=} F_1^{(x)}(s)-\Phi_{x,y}(s)\;.
    \end{array}
\end{equation}
\noindent \textbf{4. }Let us choose
\begin{equation}\label{theta}
    R>\Theta_0\stackrel{\mbox{\small def}}{=\!\!=} \frac{\mathbf{E}(\xi+\eta)^2}{2(\mathbf{E}\,\xi+\mathbf{E}\,\eta)} \left[\leqslant \frac{m_1(2)+2m_1(1)m_2(1)+m_2(2)}{2\mu_a} \right]\;;
\end{equation}
let $N$ be such that $e^{-\Lambda R}>\frac{1}{(1+NR)^{K_1\phantom{^{1}\!\!\!\!}}}\;$, and let
\begin{eqnarray*}
  q\stackrel{\mbox{\small def}}{=\!\!=}1-\left(1-\frac{\Theta_0}{R}\right)\left(e^{-\Lambda R}-\frac{1}{(1+NR)^{K_1}}\right)\varkappa(NR)\;; \nonumber\hspace{3.8cm}
  \\
  \Psi(\alpha,X_0)\stackrel{\mbox{\small def}}{=\!\!=} \left(\sum\limits_{i=0}^\infty (2i+4)^\alpha q^i\right)\left(1+ \mathbf{1}(n_0=1)2^{\alpha-1}\left(\mathrm{M}_1^{(x_0)}(\alpha)+ \mathrm{M}_2(\alpha)\right)+ \phantom{\sum\limits_{N}^N\!\!\!\!\!\!\!\!\!\!}\right. \nonumber
  \hspace{0.4cm}
  \\
   +\mathbf{1}(n_0=2)\mathrm{M}_2^{(x_0)}(\alpha)
  +2^{\alpha-1}A  \left(\frac{\alpha}{(K_1-\alpha) (K_1-\alpha-1)\mathbf{E}\,\xi}+
    \mathrm{M}_2(\alpha)\right)+
      \nonumber
      \\
  \left.\phantom{\sum\limits_{N}^N\!\!} +\frac{(1-A )\alpha}{(K_2-\alpha)(K_2-\alpha-1)\mathbf{E}\,\eta} + (i+1)\mathrm{M}_1(\alpha)+i\,\mathrm{M}_2(\alpha)\right).
\end{eqnarray*}
\section{Main result}
\textbf{Theorem 1.} {\it
Let $K>3$ and let the conditions (\ref{usl1}), (\ref{usl3}) be satisfied. Then for the process described earlier with initial state $X_0=(n_0,x_0)\;$,  for every $\alpha\in(1,K-1)$ there exists a constant $C(\alpha,X_0)<\Psi(\alpha,X_0)$ such that for all $t\geqslant 0$ the following inequality is true:
$$\big|A (t)-A \big| \leqslant\frac{C(\alpha,X_0)}{(1+t)^\alpha}\;.$$
}

\section{Proof}
\subsection{\emph{Properties of the process $X_t$}}
The process $X_t$ defined in the Subsection \ref{def} (point \textbf{2.}) is Markov.
Moreover, it possesses a strong Markov property.
We skip the standard proof of both claims.
Note that trajectories of the process $X_t$ are right continuous.
\subsection{\emph{On the stationary distribution of $X_t$}} In terms of \cite{VeZv3}, \cite{VeZv4}, the  process $X_t$ is a linear-type (piecewise linear) Markov process, and it satisfies the conditions of ergodic theorem from \cite[\S 2.6]{VeZv5} (see also \cite[Theorem 1]{VeZv6}): there exists a stationary distribution $\mathscr{P}$ on $\mathscr{X}$ such that there is a limit
$$
\lim\limits _{t\to\infty}\mathbf{P}\{n_t=j,x_t\leqslant s\}=\mathscr{P}(\{j\}\times[0,s])
$$
for any  initial state $X_0$ (again and always in the sequel $j=(1,2)$)\,\,;
\begin{equation}
 \mathscr{P}(\{j\}\times(s,\infty))= \frac{\mathbf{1}\{j=1\} \displaystyle\int\limits_s^\infty(1-F_1(s))\, \mathrm{d} s+\mathbf{1}\{j=2\} \displaystyle \int \limits_s^\infty(1-F_2(s))\,\mathrm{d} s}{\mathbf{E}\,\xi+\mathbf{E}\,\eta}\;,
\end{equation}
and $ \mathscr{P}(n_t=1)=\displaystyle\frac{\mathbf{E}\,\xi}{\mathbf{E}\,\xi+\mathbf{E}\,\eta}=A\; $.

\medskip

\subsection{\emph{Coupling method}}
To prove the Theorem  1 we will use the \emph{coupling method}, which will be now briefly recalled  (for details see \cite{VeZv7}).

Suppose some strong Markov process $X_t$  weakly converges to its (unique) stationary regime; denote its marginal distribution by  $\mathscr{P}\;$.

Suppose that on some probability space it is possible to construct two \emph{independent} versions $X_t'$ and $X_t''$ of this Markov process -- i.e., both with the same generator but possibly with different initial distributions -- such that the stopping time
$$
\tau(X_0',X_0'')\stackrel{\mbox{\small def}}{=\!\!=}\inf\{t>0:\;X_t'=X_t''\}
$$
has a finite expectation.

If, further, we have an estimate $\mathbf{E}\,\phi(\tau(X_0',X_0''))\leqslant C(X_0',X_0'')$ where $\phi(s)\uparrow$ and $\phi(s)>0$ as $s>0\;$,  then we can use a strong Markov property and \emph{coupling inequality}: for all set $\mathscr{M}\in\mathscr{B}(\mathscr{X})$
$$
       \big|\mathbf{P}\{X_t'\in \mathscr{M}\}-\mathbf{P}\{X_t''\in \mathscr{M}\}\big|\leqslant \mathbf{P}\left\{t\leqslant\tau(X_0',X_0'')\right\}=
       \mathbf{P}\left\{\phi(t)\leqslant\phi(\tau(X_0',X_0'')) \right\}\;.
$$
Hence, due to Markov's inequality,
\begin{equation}\label{kap}
       \big|\mathbf{P}\{X_t'\in \mathscr{M}\}-\mathbf{P}\{X_t''\in \mathscr{M}\}\big|\leqslant \frac{\mathbf{E}\,\phi\big(\tau(X_0',X_0'')\big)}{\phi(t)}\leqslant \frac{C(X_0',X_0'')}{\phi(t)}\;.
\end{equation}

Once the inequality (\ref{kap}) is estalished for the pair of processes, we may conclude that for the stationary process $\widetilde{X}_t$ with the initial distribution $\mathscr{P}$ and for the process $X_t$ starting from an arbitrary initial state $X_0$ we get,
\begin{equation}\label{oc1}
\begin{array}{c}
  \left|\mathbf{P}\{X_t\in \mathscr{M}\}-\mathbf{P}\{\widetilde{X}_t\in \mathscr{M}\}\right|=
   \left|\mathbf{P}\{X_t\in \mathscr{M}\}-\mathscr{P}( \mathscr{M})\}\right| \leqslant \\
   \\
  \leqslant
   \displaystyle\frac{\;\;\displaystyle\int\limits_{\mathscr{X}}C(X_0,Y) \mathscr{P}(\,\mathrm{d} Y)}{\phi(t)}=\frac{\widetilde{C}(X_0)}{\phi(t)}\;.
\end{array}
\end{equation}

Note that since the right hand side here does not depend on $\mathscr{M}\subseteq\mathscr{X}$, this inequality, of course, provides an estimate in total variation, that is,
$$
\sup_{\mathscr{M}\in \mathscr{X}}  \big|\mathbf{P}\{X_t\in \mathscr{M}\}-\mathscr{P}( \mathscr{M})\}\big| \leqslant \frac{\widetilde{C}(X_0)}{\phi(t)}\;.
$$
Also, if $ \mathscr{M}=\{n(X_t)=1\}\;$, then $\mathbf{P}\{X_t\in \mathscr{M}\}=A (t)\;$. Hence, in particular, the inequality (\ref{oc1}) implies that
$$
\big|A (t)-A \big|\leqslant\frac{\widetilde{C}(X_0)}{\phi(t)}\;.
$$

Now, the goal is to give an estimate of $\widetilde{C}(X_0)$ for the function $\phi(t)=(1+t)^\alpha\;$.

\subsection{\emph{Coupling, continued}}
We will be using a procedure first suggested in \cite{VeZv8}.
On some probability space we construct a  ``paired'' Markov process $Z_t=(Z_t',Z_t'')$ in the state space $\mathscr{X}\times\mathscr{X}$ so that the marginal distributions of the processes $Z_t'$ and $Z_t''$ coincide with the distributions of the processes $X_t'$ and $X_t''\;$, respectively:
\begin{equation}\label{Z}
    \left(Z_t',\,t\geqslant0\right) \stackrel{\mathscr{D}}{=}\left(X_t',\,t\geqslant0\right) \quad\mbox{and} \quad \left(Z_t'',\,t\geqslant0\right) \stackrel{\mathscr{D}}{=}\left(X_t'',\,t\geqslant0\right)\;;
\end{equation}
$Z_0'=X_0'$ and $Z_0''=X_0''\;$.

In addition,  we suppose, that if at  some moment $\bar{\tau}$ \, the random variable $Z_t'$ coincides with $Z_t''\;$, i.e. $Z_{\bar{\tau}}'=Z_{\bar{\tau}}''\;$, then for all $t\geqslant \bar{\tau}\;$, \, $Z_{t}'=Z_{t}''\;$.
This pair $(Z',Z'')$ is called coupling.
Of course, in general,  the processes $Z_t'$ and $Z_t''$ will be dependent.

Assuming that the process $Z_t = (Z'_t, Z''_t)$ is already constructed, let us
denote
$$
\bar{\tau}(X_0',X_0'')\left[=\bar{\tau}(Z_0',Z_0'')\right] \stackrel{\mbox{\small def}}{=\!\!=}\inf\{t>0:\;Z_t'=Z_t''\}\;.
$$
The coupling is called \emph{successful} if $\mathbf{P}\left\{\,\bar{\tau}(X_0',X_0'')<\infty\right\}=1\;$. Our coupling constructed below will be successful.

Then, we can use the coupling inequality (\ref{kap}) for the processes $Z_t'$ and $Z_t'\;$:
$$
\big|\mathbf{P}\{Z_t'\in \mathscr{M}\}-\mathbf{P}\{Z_t''\in \mathscr{M}\}\big|\leqslant \mathbf{P}\{t \leqslant\bar{\tau}(X_0',X_0'')\}\;.
$$
Due to (\ref{Z}) the same inequality holds true for $X_t'$ and $X_t''\;$.

\subsection{\emph{Construction of the process $Z_t$}}
For any distributional function $F(s)$ denote $F^{-1}(s)\stackrel{\mbox{\small def}}{=\!\!=} \inf\{u:\,F(u)=s\}$; it is well known that on the probability space $\left(\Omega^{\mathscr{L}},\mathscr{F}^\mathscr{L}, \mathbf{P}^{\mathscr{L}} \right) \stackrel{\mbox{def}}{=\!\!=}\big([0,1],\mathscr{B}([0,1]), \mathscr{L}\big)\;$ (where $\mathscr{L}\;$ is a Lebesgue measure) the random variable $\xi\stackrel{\mbox{def}}{=\!\!=}F^{-1}(\mathscr{U}_{\Omega^{\mathscr{L}}})\;$ has a distribution function $F(s)$ if $\mathscr{U}_{\Omega^{\mathscr{L}}}$ is a uniformly distributed random variable on the space $\Omega^{\mathscr{L}}\;$.

We will construct the process $Z_t\;$ on the probability space
$$
\big(\Omega,\mathscr{F},\mathbf{P}\big)
\stackrel{\mbox{def}}{=\!\!=}\prod\limits_{i=0}^\infty \Big(\big(\Omega_{i,1}^{\mathscr{L}}, \mathscr{F}_{i,1}^\mathscr{L}, \mathbf{P}_{i,1}^{\mathscr{L}} \big)\times \big(\Omega_{i,2}^{\mathscr{L}}, \mathscr{F}_{i,2}^\mathscr{L}, \mathbf{P}_{i,2}^{\mathscr{L}} \big)\Big)\;,
$$
where the probability spaces $\Omega_{i,j}^\mathscr{L}\;$ are the copies of the described above space  $\Omega^{\mathscr{L}}\;$.
The construction of $Z_t$ is based on a sequence of stopping times $t_k$, at which
$$
\mathbf{1}\left\{n(Z_{t-0}')\neq n\left(Z_{t+0}'\right)\right\}+\mathbf{1}\left\{n\left(Z_{t-0}''\right)\neq n\left(Z_{t+0}''\right)\right\}>0\;,
$$
i.e., of (random) times $t_k$ where one of the processes $Z_t'$ and $Z_t''$ -- or both of them -- changes its first component.

Let $t_0=0$ and denote
$$
m'_t\stackrel{\mbox{\small def}}{=\!\!=} n(Z'_t)\;, \quad m''_t\stackrel{\mbox{\small def}}{=\!\!=} n(Z''_t)\;, \quad z'_t\stackrel{\mbox{\small def}}{=\!\!=} x(Z'_t)\;,\quad z''_t\stackrel{\mbox{\small def}}{=\!\!=} x(Z''_t)\;.
$$
The sequence $(t_k)$ will be built by induction. Assume that $t_k$ is already determined for some $k$ and consider three cases.

\subsubsection{Case 1.} Suppose that  $Z_{t_k}'\neq Z_{t_k}''$ and $m_{t_k}'+m_{t_k}''>2$ (that is, at least one of the processes is in the set $\mathscr{S}_2$)\,.
Then on the probability space $\left(\Omega_{k,1}^\mathscr{L}\times \Omega_{k,2}^\mathscr{L}\right)\;$ we take an independent random variables $\theta_k'\stackrel{\mbox{def}}{=\!\!=}\left(F_{m_{t_k}'} ^{(z_{t_k}')}\right)^{-1} \left(\mathscr{U}_{\Omega_{k,1}^{\mathscr{L}}}\right)\;$ and \linebreak $\theta''_k\stackrel{\mbox{def}}{=\!\!=}\left(F_{m_{t_k}''} ^{(z_{t_k}'')}\right)^{-1} \left(\mathscr{U}_{\Omega_{k,2}^{\mathscr{L}}}\right)\;$, they have a distribution functions $F_{m_{t_k}'}^{(z_{t_k}')}(s)$ and $F_{m_{t_k}''}^{(z_{t_k}'')}(s)$ respectively: they are  residual times of stay of the processes $Z_t'$ and $Z_t''$ in the sets $\mathscr{S}_{m_{t_k}'}$ and $\mathscr{S}_{m_{t_k}''}\;$ correspondingly.

Denote $
\theta_k\stackrel{\mbox{\small def}}{=\!\!=} \min(\theta_k',\theta_k'')\;\mbox{ and }\;t_{k+1}\stackrel{\mbox{\small def}}{=\!\!=} t_k+\theta_k\;.
$
For $t\in[t_k,t_{k+1})$ define,
\begin{equation}\label{raspZ}
    \begin{array}{l}
      Z_t'\stackrel{\mbox{\small def}}{=\!\!=} (m_{t_k}',z_{t_k}'+t-t_k)\;; \hspace{1cm}  Z_t''\stackrel{\mbox{\small def}}{=\!\!=} (m_{t_k}'',z_{t_k}''+t-t_k)\;;
           \\ \\
      Z_{t_{k+1}}'\stackrel{\mbox{\small def}}{=\!\!=} 1\{\theta_k'=\theta_k\}\left(m'_{t_k}-(-1)^{m'_{t_k}},0\right) +
      \\
      \hspace{5cm}+1\{\theta_k'\neq\theta_k\} \left(m'_{t_k},z'_{t_k}+t_{k+1}-t_k\right)\;;
            \\ \\
      Z_{t_{k+1}}''\stackrel{\mbox{\small def}}{=\!\!=} 1\{\theta_k''=\theta_k\}\left(m''_{t_k}-(-1)^{m''_{t_k}},0\right)+
      \\
      \hspace{5cm} +1\{\theta_k''\neq\theta_k\} \left(m''_{t_k},z''_{t_k}+t_{k+1}-t_k\right) \;.
    \end{array}
\end{equation}

\subsubsection{Case 2.}
Suppose now that $Z_{t_k}'\neq Z_{t_k}''$ and $m_{t_k}'=m_{t_k}''=1\;$.
In this case, using the idea of the ``Lemma about three random variables'' (see \cite{VeZv9}) we construct on {\it one} space $\Omega_{k,1}^\mathscr{L}$ the pair of \emph{dependent} random variables $(\theta_k',\theta_k'')$ such that:
\begin{equation}\label{b}
\begin{array}{c}
  \mathbf{P}\left\{\theta_k'\leqslant s\right\}=F_1^{(z_{t_k}')}(s)\;;\quad \mathbf{P}\left\{\theta_k''\leqslant s\right\}=F_1^{(z_{t_k}'')}(s)\;;
  \\
  \\
  \mathbf{P}\{\theta_k'=\theta_k''\}= \displaystyle\int\limits_0^\infty \min\left(f_1^{(z_{t_k}')}(s), f_1^{(z_{t_k}'')}(s)'\right)\,\mathrm{d} s= \hspace{1cm}
  \\
  =\varkappa_{z_{t_k}',z_{t_k}''}\geqslant\displaystyle \int\limits_0^\infty \displaystyle\frac{K_1 e^{-\Lambda s}}{1+s+\max\left(z_{t_k}',z_{t_k}''\right)}\,\mathrm{d} s=\varkappa\left(\max\left(z_{t_k}',z_{t_k}''\right)\right)\;.
\end{array}
\end{equation}
Note that clearly $\varkappa(T)\downarrow 0\;$ if $\;T\uparrow+\infty\;$.

The construction of $(\theta_k',\theta_k'')$ is as follows.
Let
$$
\begin{array}{c}
\Xi_{x,y}(s)\stackrel{\mbox{def}}{=\!\!=}
\left\{
\begin{array}{cc}
  \Phi^{-1}_{x,y}(s) &\mbox{if }s\in[0,\varkappa_{x,y})\;;
\\
\\
\widehat{\Phi}^{-1}_{x,y}(s-\varkappa_{x,y}) &\mbox{if }s\in[\varkappa_{x,y},1)\;;
\end{array}
\right.
\\ \\
\theta_{k}'\stackrel{\mbox{def}}{=\!\!=} \Xi_{z_{\vartheta_k}',z_{\vartheta_k}''} \left(\mathscr{U}_{\Omega_{k,1}^{\mathscr{L}}}\right), \qquad \theta_{k}'\stackrel{\mbox{def}}{=\!\!=} \Xi_{z_{\vartheta_k}'',z_{\vartheta_k}'} \left(\mathscr{U}_{\Omega_{k,1}^{\mathscr{L}}}\right).
\end{array}
$$
It is easy to see that in this case the formulas (\ref{b}) for $\theta_{k}'$ and $\theta_{k}''$ are true.

Next, we again denote $t_{k+1}\stackrel{\mbox{\small def}}{=\!\!=} t_k+\min(\theta_k',\theta_k'')$ and apply the same construction given in the formulae (\ref{raspZ}). This definition and (\ref{b}) imply the inequality  $$
\mathbf{P}\left\{Z_{t_{k+1}}'=Z_{t_{k+1}}'\right\}\geqslant \varkappa\left(\max\left(z_{t_k}',z_{t_k}''\right)\right)\;.
$$

\subsubsection{Case 3.}
 Now, suppose $Z_{t_k}'=Z_{t_k}''=(m_{t_k},z_{t_k})\;$. In this case we construct random variables $\theta_k'=\theta_k''\stackrel{\mbox{def}}{=\!\!=}\left(F_{m_{t_k}}^{z_{t_k}}\right)^{-1}
 (\mathscr{U}_{\Omega_{k,1}^{\mathscr{L}}})$ (i.e., they are identical) with distribution  function $F_{m_{t_k}}^{(z_{t_k})}(s)$ on the space $\Omega_{k,1}^{\mathscr{L}}\;$, and $t_{k+1}\stackrel{\mbox{\small def}}{=\!\!=} t_k+\theta_k'\;$.  Here for $t\in[t_k,t_{k+1})\;$,
$$
 Z_t'=Z_t''=\left(m_{t_k}',z_{t_k}'+t-t_k\right)\;; \quad Z_{t_{k+1}}'=Z_{t_{k+1}}''= \left(m'_{t_k}-(-1)^{m'_{t_k}},0\right)\; .
$$

{\em This construction  gives us the desired pair $Z_t=(Z_t',Z_t'')\;$,  which satisfies (\ref{Z}) and which is suitable for the successful coupling procedure.}

Indeed, each of the processes $Z_t'$ and $Z_t''$ is an alternating, wherein periods when these processes are in the sets $\mathscr{S}_1$ or $\mathscr{S}_2$ have the distribution functions $F_1(s)$ and $F_2(s)$, respectively; the first period of their stay in the set $\mathscr{S}_{n_0'}$ or $\mathscr{S}_{n_0''}$ has a distribution function $F_{n_0'}^{(x_0')}(s)$ and $F_{n_0''}^{(x_0'')}(s)$ -- these properties are guaranteed by the construction of processes.
Moreover, for each of the processes $Z_t'\;$ and $Z_t''$  considered separately periods of its stay in the sets $\mathscr{S}_1$ and $\mathscr{S}_2$ are mutually independent.

\subsection{\emph{Using coupling method}}
Let us fix two initial values $X'_0\equiv Z_0'\neq Z_0''\equiv X''_0\;$. In this step of the proof  we will show the coupling inequality for the process $Z=(Z', Z'')\;$; hence, the same inequality will be established for the couple $(X',X'')\;$.

\subsubsection{First hits to the set $\mathscr{S}_1\;$.}
For  $t>0$ denote,
$$
\tau'(t)\stackrel{\mbox{\small def}}{=\!\!=} \inf\left\{s>t:\,Z_t'=( 1 ,0)\right\}\;, \quad \tau''(t) \stackrel{\mbox{\small def}}{=\!\!=}\inf\left\{s>t:\,Z_t''=( 1,0)\right\}\;.
$$
These are the moments of the beginning of regeneration periods for the processes $Z'$ and $Z''$ after the nonrandom $t\;$.

Denote also
$$\tau\big(Z_0',Z_0''\big)\stackrel{\mbox{\small def}}{=\!\!=} \max \big(\tau'(0), \tau''(0)\big)\;.$$

At $\tau'(0)$ the regeneration period of the process $Z'_t$ begins.

Its  length equals $\theta\stackrel{\mathscr{D}}{=}\xi+\eta$ where $\xi$ and $\eta$ were introduced in the Section 1. After that, the behaviour of $Z'$ does not depend on the initial state $Z_0'$ (given $\tau'(0)$)\,\,. The same can be said about the process~$Z_t''\;$.

Let $t>\tau\big(Z_0', Z_0''\big)\;$. Then,  there was at least one beginning of the regeneration period of each of the processes $Z'$ and $Z''$ before $t\;$.

Denote $\vartheta'(t)\stackrel{\mbox{\small def}}{=\!\!=}\left(\tau'(t)-t\right)$ -- the residual time of the last regeneration period of $Z'_t\;$, which started before time $t\;$.
From the corollary of W. Smith's Key Renewal Theorem (cf.  \cite[Theorem 2]{VeZv6}), the following inequality holds true: if  $t>\tau\left(Z_0', Z_0''\right)\;$, then
\begin{equation}\label{mo}
  \mathbf{E}\left(\vartheta'(t)\Big|t>\tau\big(Z_0', Z_0''\right)\big)\leqslant \frac{\mathbf{E}\,\theta^2}{2\mathbf{E}\,\theta}= \frac{\mathbf{E}\,(\xi+\eta)^2} {2(\mathbf{E}\,\xi+\mathbf{E}\,\eta)} \left[=\Theta_0\right]\;.
\end{equation}
The same statement applies to the process $Z''_t\;$.

Note that $\tau\big(Z_0', Z_0''\big) \leqslant \tau'(0)+\tau''(0)\;$, and, by virtue of Jensen's inequality, for all $\alpha\in(1,K-1)$
$$
\mathbf{E}\, \Big(\tau\big(Z_0', Z_0''\big)\Big)^\alpha\leqslant2^{\alpha-1} \Big(\mathbf{E}\,\big(\tau'(0)\big)^\alpha+ \mathbf{E}\,\big(\tau''(0)\big)^\alpha\Big)\;.
$$

\subsubsection{Coupling after common hit to the set  $\mathscr{S}_1\;$.}
Without loss of generality we can assume that $\tau\big(Z_0', Z_0''\big)=\tau''(0)\;$.

Let $\tau_1\stackrel{\mbox{\small def}}{=\!\!=}\tau''(0)\;$, $\tau_{k+1}\stackrel{\mbox{\small def}}{=\!\!=}\min\{\tau''(t),\,t>\tau_{k+1}\}\;$; $\{\tau_k\}$ is a sequence of beginnings of regeneration periods of $Z''\;$.

And let $\widetilde{\tau}_k\stackrel{\mbox{\small def}}{=\!\!=}\inf\{t>\tau_k:\,Z''_t=(2,0)\}$ -- time of the (first) jump of $Z_t''$ to the set $\mathscr{S}_2$ after time $\tau_k\;$.

Denote the event $\mathscr{E}_k\stackrel{\mbox{\small def}}{=\!\!=}\left\{\vartheta'(\tau_k)<R\; \&\;(\widetilde{\tau}_{k}-\tau_k)\in(R,NR)\right\}\;$, i.e. at time $\widehat{\tau}_k\stackrel{\mbox{\small def}}{=\!\!=}\tau_k+\vartheta'(\tau_k)$ \,\, $Z'_{\,\widehat{\tau}_k}=(1,0)\;$, \, $ Z''_{\,\widehat{\tau}_k}=(1,z)\;$, and $z<R\;$.

Using (\ref{mo}) and condition (\ref{usl1}) by Markov inequality we can estimate $\mathbf{P}\{\mathscr{E}_k\}\;$:
$$
    \mathbf{P}\{\mathscr{E}_k\}\geqslant\left( 1-\frac{\Theta_0}{R}\right) \left(e^{-\Lambda R}-\frac{1}{(1+NR)^{K_1}} \right)\stackrel{\mbox{\small def}}{=\!\!=}\pi(R,N)\;.
$$

Now, using (\ref{b}), we have: $\mathbf{P}\left\{Z_{\tau_{k+1}}'=Z_{\tau_{k+1}}''\right\}\geqslant \pi(R,N)\varkappa(RN)\stackrel{\mbox{\small def}}{=\!\!=} p\;$.

\subsection{\emph{Completion of the proof}}
 The number of regeneration periods of $Z_t''$ before the processes  $Z_t'$ and $Z_t''$ meet each other according to the scheme from the step \textbf{4.6.} (that is, any meeting outside this scheme is ignored) is a random variable $\nu$ dominated by another one with a geometric distribution with parameter $p$ ($\nu$ itself has a more complicated distribution).
Denote $q\stackrel{\mbox{\small def}}{=\!\!=}1-p\;$, and $\varsigma\left(Z_0',Z_0''\right)\stackrel{\mbox{\small def}}{=\!\!=} \inf\left\{t>0:\,Z_t'=Z_t''\right\}\;$.
Obviously,  $\varsigma\left(Z_0',Z_0''\right)\leqslant \tau_\nu\;$.

Since we know the distribution of $\tau=\tau\left(Z_0',Z_0''\right)$ and $\theta\stackrel{\mathscr{D}}{=}\xi+\eta\;$, we can obtain an estimation of $\mathbf{E}\left(1+\varsigma\left(Z_0',Z_0''\right)\right)^\alpha\;$ for all $\;\alpha\in(1,K-1)\;$: by Jensen's inequality we get,
\begin{eqnarray}
\mathbf{E} \Big(1+\varsigma\left(Z_0', Z_0''\right)\Big)^\alpha\leqslant
\hspace{8.5cm}
\nonumber
\\
\leqslant\mathbf{E}\left(1+ \tau\left(Z_0', Z_0''\right)+\xi+
\sum\limits_{i=1}^\infty \left(\mathbf{P}\{\nu=i\}\sum\limits_{k=1}^i     (\xi_k+\eta_k)\right)\right)^\alpha\leqslant
\hspace{2cm}
\nonumber
\\
\leqslant \sum\limits_{i=1}^\infty q^{i-1}\mathbf{E}\left(1+\tau'(0)+ \tau''(0)+\xi+ \sum\limits_{k=1}^i(\xi_k+\eta_k)\right)^\alpha     \leqslant
\hspace{2.5cm}
\nonumber
\\
\leqslant \sum\limits_{i=1}^\infty q^{i-1}(2i+4)^{\alpha-1} \bigg(1+\mathbf{E}\big(\tau'(0)\big)^\alpha +
\hspace{2.5cm}
\nonumber
\\
+\mathbf{E}\big(\tau''(0)\big)^\alpha+ (i+1)\mathbf{E}\,\xi^\alpha+ i\,\mathbf{E}\,\eta^\alpha \bigg)\leqslant
\hspace{2cm}
\nonumber
\\
\leqslant \sum\limits_{i=1}^\infty q^{i-1}(2i+4)^{\alpha-1}     \Bigg(1+
\hspace{6.5cm}
\nonumber
\\
+\mathbf{1}\left(n_0'=1\right)2^{\alpha-1} \Big(\mathrm{M}_1^{\left(x_0'\right)}(\alpha)+     \mathrm{M}_2(\alpha)\Big)+\mathbf{1}\left(n_0'=2\right) \mathrm{M}_2^{\left(x_0'\right)}(\alpha)+
\hspace{1cm}
\nonumber
\\
+ \mathbf{1}\left(n_0''=1\right)2^{\alpha-1} \Big(\mathrm{M}_1^{(x_0'')}(\alpha)+ \mathrm{M}_2(\alpha)\Big)+     \mathbf{1}\left(n_0''=2\right) \mathrm{M}_2^{\left(x_0''\right)}(\alpha)+
\hspace{1cm}
\nonumber
\\
+ (i+1)\mathrm{M}_1(\alpha)+i\mathrm{M}_2(\alpha)\Bigg)=
\hspace{1cm}
\nonumber
\\
=C\left(\alpha,Z_0',Z_0''\right)= C\left(\alpha,X_0',X_0''\right)\;.
\hspace{3cm}
\label{15}
\end{eqnarray}

Now, considering (\ref{oc1}), it is necessary to estimate the integral \linebreak $\displaystyle\int\limits_{\mathscr{X}} C\left(\alpha,X_0',X_0''\right)\, \mathscr{P}\left(\mathrm{d}\,X_0''\right) \;$.

First, we make an estimations of the integrals:
\begin{eqnarray}
\int\limits_{\mathscr{X}} \left( \mathbf{1}\left(n_0''=1\right)\right) \mathscr{P}\left(\mathrm{d}\,X_0''\right)= \mathscr{P}(\{1\}\times(0,+\infty))=A\;;
\hspace{3cm}
\label{1}
\\
\int\limits_{\mathscr{X}} \left( \mathbf{1}\left(n_0''=2\right)\right) \mathscr{P}\left(\mathrm{d}\,X_0''\right)= \mathscr{P}(\{1\}\times(0,+\infty))=1-A\;;
\hspace{2.5cm}
\label{2}
\\
\int\limits_{\mathscr{X}} \left( \mathbf{1}\left(n_0''=1\right) \mathrm{M}_1^{(x_0'')}(\alpha)\right)
\mathscr{P}\left(\mathrm{d}\,X_0''\right)=
\nonumber
\hspace{5.3cm}
\\
=\int\limits_{\mathscr{X}} \left( \frac{\mathbf{1}\left(n_0''=1\right)\alpha} {1-F_1(x_0'')}\int\limits_0^\infty \frac{s^{\alpha-1}}{(1+s+x_0'')^{K_1{\phantom{^1}}}\!\!\!} \mathrm{d}\, s\right) \times
\nonumber
\hspace{3cm}
\\
\times\mathrm{d}\left(1-\frac{\int\limits_{x_0''}^\infty (1-F_1(s))\,\mathrm{d}\,s}{\mathbf{E}\,\xi+\mathbf{E}\,\eta} \right)<
\nonumber
\hspace{2.5cm}
\\
<
\int\limits_{0}^\infty\left(\frac{\alpha}{1-F_1(x_0'')}\int\limits_0^\infty \frac{(1+s+x_0)^{\alpha-1}}{(1+s+x_0'')^{K_1{\phantom{^1}}}\!\!\!} \,\mathrm{d}\, s\right) \frac{1-F_1(x_0'')}{\mathbf{E}\,\xi+\mathbf{E}\,\eta} \,\mathrm{d}\, x_0''\leqslant
\nonumber
\\
\leqslant
\int\limits_{0}^\infty \frac{\alpha(1+x_0'')^{\alpha-K_1}} {(K_1-\alpha) (\mathbf{E}\,\xi+\mathbf{E}\,\eta)}\,\mathrm{d}\,x_0''= \frac{\alpha}{(K_1-\alpha)(K_1-\alpha-1) (\mathbf{E}\,\xi+\mathbf{E}\,\eta)}=
\nonumber
\\
=\frac{\alpha A}{(K_1-\alpha)(K_1-\alpha-1) \mathbf{E}\,\xi}\;;
\hspace{2cm}
\label{3}
\\
\mbox{analogously }\int\limits_{\mathscr{X}} \left( \mathbf{1}\left(n_0''=2\right) \mathrm{M}_1^{(x_0'')}(\alpha)\right)
\mathscr{P}\left(\mathrm{d}\,X_0''\right)<
\nonumber
\hspace{3.5cm}
\\
<\frac{\alpha(1-A)}{(K_2-\alpha)(K_2-\alpha-1) \mathbf{E}\,\eta}\;.
\hspace{2cm}
\label{4}
\end{eqnarray}

Considering, that for all constant $\Upsilon$ and set $\mathscr{M}\in\mathscr{B}(\mathscr{X})$ we have \linebreak $\displaystyle\int\limits_{\mathscr{M}} \Upsilon\, \mathscr{P}\left(\mathrm{d}\,X_0''\right)= \Upsilon\mathscr{P}(\mathscr{M})\;$, and collecting formulas (\ref{1})--(\ref{4}), we can estimate the integral $\displaystyle\int\limits_{\mathscr{X}} C\left(\alpha,X_0',X_0''\right)\mathscr{P}\left(\,\mathrm{d} X_0''\right)\;$, end this estimation give as an inequality
\begin{equation}\label{sovs}
    \int\limits_{\mathscr{X}} C\left(\alpha,X_0',X_0''\right)\mathscr{P}\left(\,\mathrm{d} X_0''\right) \leqslant \Psi\left(\alpha,X_0'\right)\;,
\end{equation}
which completes the proof of the theorem.

\textbf{Remark.} The estimate (\ref{sovs}) could be improved; moreover, a more careful choice of  parameters $R$ and $N$ may provide some enhancement of this bound.

{\bf Acknowledgments.} The authors are grateful to L.~G.~Afanasieva and V.~V.~Kozlov for very useful consultations.

\end{document}